\documentclass[journal,a4paper]{IEEEtran}
\usepackage[T1]{fontenc}% optional T1 font encoding
\usepackage{layouts}
\usepackage{cite}
\usepackage{algorithmic}
\usepackage{array}
\ifCLASSOPTIONcompsoc
  \usepackage[caption=false,font=normalsize,labelfont=sf,textfont=sf]{subfig}
\else
  \usepackage[caption=false,font=footnotesize]{subfig}
\fi

\usepackage{float}
\usepackage{verbatim}
\usepackage{bm}% bold math
\usepackage{amsmath,amsfonts,amsthm,amssymb,times,mathrsfs,txfonts}
\usepackage{enumitem}
\usepackage{multirow}
\usepackage[pdftex]{graphicx}
\DeclareGraphicsExtensions{.pdf,.jpeg,.png}

\usepackage[pdftex,colorlinks=true,bookmarks=true,citecolor=blue,urlcolor=blue]{hyperref} %pdflatex

\newcommand{\field}[1]{\mathbb{#1}}
\newcommand{\fs}[1]{\mathsf{#1}}

\newcommand{\tp}{\intercal}% transpose operation

\newcommand{\ovl}[1]{\overline{#1}}
\newcommand{\bigO}[1]{\mathop{O}\left(#1\right)}

\newcommand{\vv}[1]{\boldsymbol{#1}}
\newcommand{\vs}[1]{\boldsymbol{#1}}

\DeclareMathOperator{\Si}{Si}
\DeclareMathOperator{\Cin}{Cin}

\DeclareMathOperator{\sinc}{sinc}

\newtheorem{theorem}{Theorem}[section]
\newtheorem{prop}[theorem]{Proposition}

\newtheorem{lemma}[theorem]{Lemma}

\newtheorem{rem}{Remark}[section]

\newcommand{\et}{\textit{et~al.}}

\begin{document}

\title{New method of bandlimited extrapolation}

\author{Vishal Vaibhav% <-this % stops a space
\thanks{Vishal Vaibhav is with the Department of Physics, 
Indian Institute of Technology Delhi, 
Hauz Khas, New Delhi 110016, India, Email:~\texttt{vvaibhav@iitd.ac.in}}
}
\maketitle

\begin{abstract}
The paper deals with numerical solution of the Fredholm integral equation 
associated with the classical problem of extrapolating
bandlimited functions known on $(-1,1)$ to the entire real line. The approach
presented can be characterized as the degenerate kernel method using the 
spherical Bessel functions as basis functions where the Tikhonov regularization
is applied at the discrete level in order to deal with the ill-posedness of the
problem. 
\end{abstract}

\section{Introduction}
%%%%%%%%%%%%%%%%%%%%%%%%%%%%%%%%%%%%%%%%%%%%%%%%%%%%%%%
In this paper, we revisit the classical problem of
extrapolation of bandlimited signals known on $(-1,1)$ to the entire
real line. The earliest known method of solving this problem is the celebrated 
Papoulis-Gerchberg (PG) algorithm~\cite{G1974,Papoulis1975} which takes an iterative
approach. Since the publication of these seminal works, a number of new
techniques have been developed to overcome the shortcomings of the original PG
algorithm~\cite{SS1978,Cadzow1979,XC1983,XC1984,Chen2006}. Our approach here is 
based on the discretization of the 
associated Fredholm integral equation in 
the spirit of the work of Khare and George~\cite{KG2003}, and, Walter and 
Solesky~\cite{WS2005}. The key difference here lies in 
the fact that instead of translates of the sinc function (i.e. the 
Whittaker-Kotelnikov-Shannon sampling basis functions), we propose to use 
the spherical Bessel functions for the discretization of the integral equation. This 
representation is then used to solve the classical problem of
extrapolation of bandlimited signals known on $(-1,1)$ to the entire real 
line. It is well known that bandlimited extrapolation 
is inherently an ill-posed problem; therefore, we must turn to the standard techniques 
of solving ill-posed problems such as 
Tikhonov regularization~\cite{TA1977,H1998,N1998}. In contrast
to the approach adopted in~\cite{Chen2006}, the regularization is applied to the
discrete system. It must be emphasized here that in the work~\cite{KG2003,WS2005}, the
authors propose the method as a way of computing the eigenfunctions which are
known to be the prolate spheroidal wave functions of order zero (PSWFs). As noted 
in~\cite{Kirby2006}, the discrete system obtained via the associated Fredholm 
equation is ill-conditioned and provides poor accuracy compared to the standard 
method of computing PSWFs on $(-1,1)$ due to Bouwkamp~\cite{B1947} which seeks 
to solve a related Sturm-Liouville problem 
(see~\cite{XRY2001,ORX2013} for an extensive treatment of the method). Therefore, 
the significance of the work~\cite{KG2003,WS2005} must be viewed in the light
of the fact that they provide a discrete framework for the solution of the 
inverse problem considered in the paper. At the heart of the method described 
in~\cite{KG2003,WS2005} and in this paper is the fact that 
degenerate approximation of the sinc-kernel 
(either in terms of the translates of sinc function as in~\cite{KG2003,WS2005} or 
the spherical Bessel functions proposed in this paper) can be obtained without 
any computational effort. 

Among the more recent attempts to address the bandlimited extrapolation problem,
we refer the reader to~\cite{G2010} where the author uses PSWFs for smaller 
frequencies and a compressed sampling approach for the higher frequencies of 
the signal. Let us briefly remark that the spherical Bessel functions would 
appear naturally here in the low-frequency case had the author used the exact
integral as opposed to the Gaussian quadrature. Note that for highly oscillatory
integrals, Gaussian quadrature is not the most preferred method of numerical 
integration. For the sake of brevity of presentation, we do not attempt to provide 
any further review of the existing literature on duration and bandwidth 
limiting; the reader may find an exhaustive survey with applications 
in Hogan and Lakey~\cite{HL2012}.

As far as the numerical technique is concerned, the method presented in this 
paper can be classified as the degenerate
kernel method which is one of the standard methods for solving Fredholm
equations~\cite{A2009}. Among the other methods are the Nystr\"om
method and the projection method (see~\cite{A2009} for a comprehensive
treatment); however, our method has an intuitive 
appeal on account of the fact that the solutions are bandlimited functions, 
and, it is natural to think of their representation in terms of bandlimited 
functions that form an orthonormal basis. In particular, our objective is to
compare the performance of the spherical Bessel functions with that of 
translates of sinc function. The tests reveal that spherical Bessel functions
exhibit better performance for the class of functions in $\fs{B}_{\sigma}^2$
(Bernstein spaces~\cite{N1975}) whose spectrum belongs to $\fs{C}^{p}(-\sigma,\sigma)$ 
where the integer $p\geq1$. Finally, let us mention that, in treating the
extrapolation problem, we have assumed that the input to the algorithm is not 
contaminated with noise. Note that the inverse problem in question is ill-posed 
regardless of the presence of noise and it is of interest to understand the 
fundamental limitations of the numerical algorithm in treating such 
problems. Therefore, we assume that it is possible to
make precise measurements of the input or make a sufficiently accurate 
estimate of the input through multiple measurements. The input signal 
can then be obtained as a fitted polynomial over the domain $(-1,1)$. The possibility of 
treating the noisy case in the discrete framework considered in this paper 
is deferred to future research.

The original motivation 
behind this work was to facilitate the estimation of compactly supported 
scattering potentials in inverse scattering problems when a complete set of 
scattering data is unavailable (see~\cite{V2018TL} for the characterization 
of the scattering coefficients when the scattering potential is compactly 
supported). In particular, the problem of estimation of 
the fiber Bragg gratings~\cite{FZM1999,SW2003} is one such 
area of application where the ideas being presented in this paper are directly 
applicable. A similar application would be in estimation of the refractive 
index profile of a layered media through inverse 
scattering~\cite{Y1989,CR1992}. Further, the ideas presented in this paper 
can be readily extended to treat inverse problems of general nature 
where the kernel of the Fredholm equation is bandlimited~\cite{KG2005,K2007}.

The paper is organized as follows: Sec.~\ref{sec:PSWFs} discusses the 
discretization of the Fredholm equation associated with the bandlimited
extrapolation problem which is considered in Sec.~\ref{sec:Xtrap-BL}. 
Sec~\ref{sec:conclusion} concludes the paper.

\section{Angular Prolate Spheroidal Wave Functions}
\label{sec:PSWFs}
Let $\Omega=(-1,1)$. For a given $\sigma>0$, referred to as the
\emph{bandlimiting} parameter, the prolate spheroidal wave functions 
(PSWFs)~\cite{SPI1961,S1983} are defined as eigenfunctions of the 
eigenvalue problem given by
\begin{equation}\label{eq:PSWFs_EVP1}
\int_{\Omega}\frac{\sin[\sigma(t-s)]}{\pi(t-s)}\phi(s)ds=\lambda\phi(t),\quad
t\in\Omega.
\end{equation}
The equation~\eqref{eq:PSWFs_EVP1} is also valid for 
all $t\in\field{R}$ and it defines the values of the PSWFs on the real
line. The eigenvalues are all positive real numbers indexed in the decreasing
order of their magnitude, $1>\lambda_0>\lambda_1>\ldots>\lambda_n>\ldots,$
such that $\lim_{n\to\infty}\lambda_n=0$ and the corresponding eigenfunctions are denoted by 
$\phi_n(t),\,\,n=0,1,\ldots$, respectively. The eigenfunctions form a complete orthogonal 
basis in $\fs{L}^2(\Omega)$. On the real line, they form a complete orthonormal basis spanning
the class of bandlimited functions in $\fs{L}^2(\field{R})$. The double
orthogonality property of PSWFs is characterized by the following relations:
\begin{equation}\label{eq:inner-prod}
\begin{split}
&\langle\phi_n,\phi_m\rangle_{\Omega}=\int_{\Omega}\phi_n(t)\phi_m(t)dt=\lambda_n\delta_{mn},\\
&\langle\phi_n,\phi_m\rangle_{\field{R}}=\int_{\field{R}}\phi_n(t)\phi_m(t)dt=\delta_{mn},
\end{split}
\end{equation}
where $\langle\cdot,\cdot\rangle$ stands for the inner product and $\delta_{mn}$
denotes the Kronecker delta. We enumerate some of the important properties of 
PSWFs for ready reference~\cite{ORX2013}: The eigenfunctions $\phi_n(t)$ satisfy the parity relation
$\phi_n(-t) = (-1)^n\phi_n(t)$ and have exactly $n$ zeros in $\Omega$. The 
PSWFs also satisfy the eigenvalue problem
\begin{equation}\label{eq:PSWFs_EVP2}
\int_{\Omega}e^{i\sigma st}\phi_n(s)ds=\nu_n\phi_n(t),\quad t\in\Omega,
\end{equation}
where $\nu_n=i^n\sqrt{2\pi\lambda_n/\sigma}$. The equation~\eqref{eq:PSWFs_EVP2}
is also valid for $t\in\field{R}$ so that either of the 
relationships~\eqref{eq:PSWFs_EVP1} or~\eqref{eq:PSWFs_EVP2} can be
used for computing the values of the PSWFs outside $\Omega$. The Fourier transform of 
$\phi_n(t)$, denoted by $\Phi_n(\xi)$, is given by
\begin{equation}
\Phi_n(\xi)
=\int_{\field{R}}\phi_n(t)e^{-i\xi t}dt\\
=(-i)^n\sqrt{\frac{2\pi}{\sigma\lambda_n}}\phi_n\left(\frac{\xi}{\sigma}\right)
\Pi\left(\frac{\xi}{\sigma}\right),
\end{equation}
where $\Pi(\xi)$ is the rectangle function which is defined to be unity for
$\xi\in\Omega$ and zero otherwise.

\begin{rem}\label{rem:I}
In most of the extrapolation problems, the value of $\phi_n(t)$ is needed in the
region $t\in\field{R}\setminus(-1,1)$. One of the standard methods to compute
these values is to apply numerical quadrature techniques to the oscillatory
integral in~\eqref{eq:PSWFs_EVP2} which reads as
\begin{equation}
\phi_n(t)=\frac{(-i)^n}{\sqrt{2\pi\lambda_n/\sigma}}
\int_{-1}^1\phi_n(s)e^{i\sigma st}ds.
\end{equation}
Gauss-type quadrature schemes tend to perform poorly in computing these 
integrals on account of the oscillatory nature of the integrand which
deviates considerably from polynomials, specially for larger values of
$\sigma$. There is a vast amount of literature devoted to treating such
problems, for instance, see~\cite[Section~2.10]{DR1984}) and the references
therein. 
\end{rem}
The problem of determination of $\phi_n(t)$ for $t\in\field{R}\setminus(-1,1)$
requires us to deal with the oscillatory integral in~\eqref{eq:PSWFs_EVP2}. The most 
convenient method that fits naturally into our setting 
is the method of Bakhvalov and Vasil'eva~\cite{BV1968} which, employing the series
expansion in the Bouwkamp method~\cite{B1947} in terms of normalized Legendre
polynomials ($\ovl{P}_n(t)=\sqrt{n+1/2}P_n(t)$),
\begin{equation}\label{eq:LegendreSeries}
\phi_n(t)= \sum_{m=0}^{\infty}\alpha^{(n)}_m \ovl{P}_m(t),
\end{equation}
yields the representation
\begin{equation}\label{eq:BesselSeries}
\phi_n(t)=\sqrt{\frac{2\sigma}{\pi\lambda_n}}\sum_{m=0}^{\infty}
i^{m-n}\alpha^{(n)}_m\sqrt{m+{1}/{2}}\,j_{m}(\sigma t),
\end{equation}
using the identity
\begin{equation}
\int_{-1}^{1}e^{ist}P_n(s)ds=2i^{n}j_n(t),
\end{equation}
where $j_n(t)$ denotes the spherical Bessel function of the first kind defined
as $j_n(t) = \sqrt{\pi/2t}J_{n+1/2}(t)$. On account of the presence of the 
factor $\lambda_n^{-1/2}$, the numerical conditioning of the
expression~\eqref{eq:BesselSeries} is extremely poor for 
$n>2\sigma/\pi$~\cite{XRY2001}. Therefore, we would like to avoid computing the 
PSWFs and develop a representation directly in terms of the spherical Bessel
function. Note that this does not remedy the ill-conditioning issues
because it stems from the inherently ill-posed nature of the problem.

In the remainder of this section, we consider two ways of discretizing the
Fredholm integral equation~\eqref{eq:PSWFs_EVP1} with the objective of designing
a regularization approach within the discrete system in order to solve the
classical bandlimited extrapolation problem. 
%%%%%%%%%%%%%%%%%%%%%%%%%%%%%%%%%%%%%%%%%%%%%%%%%%%%%%%%%
%%%%%%%%%%%%%%%%%%%%%%%%%%%%%%%%%%%%%%%%%%%%%%%%%%%%%%%
\subsection{The approach based on spherical Bessel functions}
\label{sec:PSWFs-Bessel-IE}
Let us introduce the normalized spherical Bessel functions
\begin{equation}
\ovl{j}_n(t)=\sqrt{\frac{\sigma(2n+1)}{\pi}}j_n(\sigma t),
\end{equation}
for convenience so that the orthonormality relation reads as 
$\langle\ovl{j}_m,\ovl{j}_n\rangle_{\field{R}} =\delta_{mn}$ 
for all integers $m,n\geq0$. The method outlined in this section is motivated 
by the following representation of the sinc-kernel in~\eqref{eq:PSWFs_EVP1}:
\begin{equation}\label{eq:sinc_sphj}
\frac{\sin{\sigma(t-s)}}{\pi(t-s)} 
=\sum_{n=0}^{\infty}\ovl{j}_n(t)\ovl{j}_n(\sigma s).
\end{equation}
The eigenvalue problem in~\eqref{eq:PSWFs_EVP1} can be 
discretized by writing
\begin{equation}\label{eq:psf_bessel}
\phi(t)=\sum_{n=0}^{\infty}\beta_n\ovl{j}_n(t),    
\end{equation}
where $\beta_n$ are unknown coefficients and using~\eqref{eq:sinc_sphj} so that
\begin{equation*}
\lambda\sum_{m=0}^{\infty}\beta_m\ovl{j}_m(t)
=\sum_{m=0}^{\infty}\ovl{j}_m(t)\sum_{n=0}^{\infty}
\left(\int_{-1}^{1}\ovl{j}_m(s)\ovl{j}_n(s)\right)\beta_n.    
\end{equation*}
Equating the coefficient of $\ovl{j}_m(t)$ on both sides, we have
\begin{equation}\label{eq:symm_evp}
\lambda\beta_m 
= \sum_{n=0}^{\infty}J_{mn}\beta_n,   
\end{equation}
where the matrix elements $J_{mn}$ are given by
\begin{equation}\label{eq:jmn}
J_{mn}=\langle\ovl{j}_m,\ovl{j}_n\rangle_{\Omega}
=\int_{-1}^{1}ds\,\ovl{j}_m(s)\ovl{j}_n(s).    
\end{equation}
For finite $\sigma$, using Cauchy-Schwarz inequality, we have 
$|J_{mn}|\leq\|\ovl{j}_m\|_{\fs{L}^2(\Omega)}\|\ovl{j}_n\|_{\fs{L}^2(\Omega)}<1$
where we used the fact that
$\|\ovl{j}_m\|_{\fs{L}^2(\Omega)}<\|\ovl{j}_m\|_{\fs{L}^2(\field{R})}=1$. It must be 
noted that $J_{mn}=0$ for $m \neq n\mod{2}$; therefore, the even system
of equations can be decoupled with that of the odd. For $m\neq n$ and $m+n=0\mod{2}$, we have
\begin{multline}\label{eq:jm-neq-n}
J_{mn}=\frac{2\sigma^2\sqrt{(2m+1)(2n+1)}}{\pi(m+n+1)}\\
\times\left[\left(\frac{j_{m-1}(\sigma)j_n(\sigma)
-j_{n-1}(\sigma)j_m(\sigma)}{m-n}\right)-\frac{1}{\sigma}j_{m}(\sigma)j_n(\sigma)\right].
\end{multline}
However, the diagonal elements $J_{mm}$ cannot be computed in a closed form. 
Let us now state the following decay property of the overlap integrals:
\begin{lemma}
The overlap integrals defined by~\eqref{eq:jmn} satisfies the estimate
\begin{equation}\label{eq:jmn-estimate2}
|J_{mn}|\leq \frac{e\sigma}{2}\left(\frac{e\sigma}{2m+2}\right)^m
\left(\frac{e\sigma}{2n+2}\right)^n.
\end{equation}
\end{lemma}
\begin{proof}
To prove this lemma, we start with the integral representation
\begin{equation}
j_n(x) = \frac{x^n}{2^{n+1}n!}\int_0^{\pi}\cos(x\cos\theta)\sin^{2n+1}\theta d\theta,
\end{equation}
so that
\begin{multline*}
J_{mn}=C_{mn}\int_0^{\pi}\int_0^{\pi}d\theta'
d\theta\int_{-1}^{1}ds\,{s^{n+m}}\\
\times\cos(\sigma s\cos\theta)
\cos(\sigma s\cos\theta')
\sin^{2m+1}\theta\sin^{2n+1}\theta',
\end{multline*}
where 
\begin{equation*}
C_{mn}=\frac{\sigma^{m+n+1}\sqrt{(2m+1)(2n+1)}}{2^{n+m+2}n!m!\pi}.
\end{equation*}
This leads to the following estimate
\begin{equation}\label{eq:jmn-estimate1}
\begin{split}
|J_{mn}|
&\leq\frac{\sigma^{m+n+1}\sqrt{(2m+1)(2n+1)}}{2^{n+m+1}n!m!(m+n+1)\pi}
\frac{2^{m+1}m!}{(2m+1)!!}\frac{2^{n+1}n!}{(2n+1)!!}\\
&=\frac{\sqrt{(2m+1)(2n+1)}}{\pi(m+n+1)}\frac{2\sigma^{m+n+1}}{(2m+1)!!(2n+1)!!},
\end{split}
\end{equation}
where we have used the result 
\begin{equation}
\int_0^{\pi}(\sin\theta)^{2n+1} d\theta = \frac{2^{n+1}n!}{(2n+1)!!}.
\end{equation}
Using the inequality~\cite[Thm.~1.6]{B2008}
\begin{equation}
\frac{1}{(2n+1)!!}=\frac{\sqrt{\pi}}{2^{n+1}\Gamma(n+3/2)}
\leq\sqrt{\frac{\pi}{2e}}\left(\frac{e}{2n+2}\right)^{n+1},
\end{equation}
to~\eqref{eq:jmn-estimate1} yields the estimate 
in~\eqref{eq:jmn-estimate2}.

\end{proof}

This inequality allows us to conclude that $J$ is a Hilbert-Schmidt operator on
$\ell^2$:
\[
\|J\|^2_{\text{F}}=\sum_{m=0}^{\infty}\sum_{m=0}^{\infty}|J_{mn}|^2\leq
\left(\frac{e\sigma}{2}\right)^2
\left[\sum_{n=0}^{\infty}\left(\frac{e\sigma}{2n+2}\right)^{2n}\right]^2,
\]
where `F' stands for the Frobenius norm. It is also possible to estimate the
spectral norm of $J$ as follows:
\begin{align*}
\sum_{m,n\geq0}J_{mn}\beta_m\beta_n
&=\int^1_{-1}\left(\sum_{n\geq0}\beta_n\ovl{j}_n(s)\right)^2ds\\
&<\int^{\infty}_{-\infty}\left(\sum_{n\geq0}\beta_n\ovl{j}_n(s)\right)^2ds
=\|\vs{\beta}\|^2_{\ell^2},
\end{align*}
so that $\|J\|_s<1$ provided $\sigma<\infty$. These properties are summarized 
in the following proposition:
\begin{prop}
The infinite matrix $J$ defines a self-adjoint, positive definite and 
compact operator on $\ell^2$.
\end{prop}
Thus, $J$ admits of an orthonormal sequence of eigenvectors which 
correspond to positive eigenvalues less than unity. Given the real nature of the 
matrix $J$, we may further assume without loss of generality that its eigenvectors 
are real. Let $\vs{\beta}_n=(\beta^{(n)}_0,\beta^{(n)}_1,\ldots)^{\tp}$ define 
$\phi_n(t)$ in the sense of~\eqref{eq:psf_bessel}. We assume that the eigenvectors 
$\vs{\beta}_n$ are normalized such that
$\|\vs{\beta}_n\|_{\ell^2}=1$. The inner-product on $\field{R}$ is given by
$\langle\phi_m,\phi_n\rangle_{\field{R}}=\vs{\beta}^{\tp}_m\vs{\beta}_n=\delta_{mn}$
using the orthogonality of the spherical Bessel functions. The inner-product 
on $\Omega$ is given by
$\langle\phi_m,\phi_n\rangle_{\Omega}=\vs{\beta}^{\tp}_mJ\vs{\beta}_n = \lambda_n\delta_{mn}$.
The eigenfunctions $\phi_n(t)$ have definite parity which implies that 
the Fourier transform of $\phi_n(t)$ must be either purely real or
imaginary. Using the identity
\begin{equation}
\int_{\field{R}}dt\,\ovl{j}_n(t)e^{-i\xi t}=\sqrt{\frac{2\pi}{\sigma}}
(-i)^n\ovl{P}_n\left(\frac{\xi}{\sigma}\right),\quad\xi\in(-\sigma,\sigma),
\end{equation}
we have
\begin{equation}
{\Phi}_n(\xi)
=\sqrt{\frac{2\pi}{\sigma}}\sum_{k=0}^{\infty}\beta^{(n)}_k(-i)^k
\ovl{P}_k\left(\frac{\xi}{\sigma}\right)
\Pi\left(\frac{\xi}{\sigma}\right).
\end{equation}
Noting that, for non-zero value of the coefficients, $k$ either runs 
through all even or all odd values; 
it is evident that $\beta^{(n)}_k\in\field{R}$ (or identically zero) which confirms that the 
eigenfunctions are real valued. 

\begin{lemma}
The function ${\Phi}_n(\sigma s),\,s\in\Omega,$ also satisfies 
the eigenvalue problem in~\eqref{eq:PSWFs_EVP1}. 
\end{lemma}
\begin{proof}
Consider
\begin{equation}\label{eq:Phi}
\begin{split}
\mathcal{I}&=\int^{1}_{-1}\frac{\sin[\sigma(t-s)]}{\pi(t-s)}{\Phi}_n(\sigma s)ds\\
&=\sigma\int^{1}_{-1}\int_{-1}^{1}e^{i\sigma\eta(s-t)}\frac{d\eta}{2\pi}{\Phi}_n(\sigma s)ds\\
&=\sum_{k=0}^{\infty}\beta^{(n)}_k\int_{-1}^{1}d\eta e^{-i\sigma\eta t} \ovl{j}_k(\eta).
\end{split}
\end{equation}
Observing that~\cite{Olver:2010:NHMF}
\begin{equation}\label{eq:exp_jn}
e^{-i\sigma\eta t}=\sqrt{\frac{2\pi}{\sigma}}
\sum_{l=0}^{\infty}(-i)^l\ovl{j}_l(\eta)\ovl{P}_l(t),\quad t\in\Omega,
\end{equation}
the proof follows from
\begin{align*}
\mathcal{I}&=\sqrt{\frac{2\pi}{\sigma}}\sum_{k=0}^{\infty}\beta^{(n)}_k
\sum_{l=0}^{\infty}(-i)^l\ovl{P}_l(t)
\int_{-1}^{1}d\eta\,\ovl{j}_k(\eta)\ovl{j}_l(\eta)\\
&=\sqrt{\frac{2\pi}{\sigma}}\sum_{l=0}^{\infty}\sum_{k=0}^{\infty}\beta^{(n)}_kJ_{kl}(-i)^l\ovl{P}_l(t)\\
&=\lambda_n\sqrt{\frac{2\pi}{\sigma}}\sum_{l=0}^{\infty}\beta^{(n)}_l(-i)^l\ovl{P}_l(t)
=\lambda_n{\Phi}_n(\sigma t),\quad t\in\Omega. 
\end{align*}
\end{proof}
This lemma shows that we can write $\Phi_n(\sigma t)=\gamma_n\phi_n(t)$ on account 
of the fact that the eigenfunctions are unique up to scalar multiplier. Noting that
\begin{align*}
&\gamma^{-2}_n\int_{-1}^1{\Phi}^2_n(\sigma t)dt=
\gamma^{-2}_n(2\pi/\sigma)\sum_k\beta^{(n)}_k\beta^{(n)}_k(-i)^{2k}=\lambda_n,
\end{align*}
it follows that for the even case $\gamma_n=\pm\sqrt{2\pi/\sigma\lambda_n}$ and for the odd case
$\gamma_n=\pm i\sqrt{2\pi/\sigma\lambda_n}$. The ambiguity in sign is resolved by
requiring $\phi_n(t)\rightarrow P_n(t)$ as $\sigma\rightarrow0$. This leads to
$\gamma_n=(-i)^n\sqrt{2\pi/\sigma\lambda_n}$. From~\eqref{eq:Phi} it
also becomes evident that 
\begin{equation*}
(-1)^n\lambda_n\Phi_n(\sigma t)
=\int_{-1}^{1}d\eta e^{i\sigma\eta t}
\sum_{k=0}^{\infty}\beta^{(n)}_k\ovl{j}_k(\eta)
=\int_{-1}^{1}d\eta e^{i\sigma\eta t}\phi_n(\eta).
\end{equation*}
The eigenvalues $\nu_n$ works out to be
$\nu_n=(-1)^n\lambda_n\gamma_n=i^n\sqrt{2\pi\lambda_n/\sigma}$. This proves the following lemma
\begin{lemma}
The functions $\phi_n(t)$ and, equivalently, $\Phi_n(\sigma t)$ satisfy the
eigenvalue problem in~\eqref{eq:PSWFs_EVP2} for $t\in\Omega$. 
\end{lemma}

\subsection{The approach based on sampling theory}
\label{sec:PSWFs-sinc-IE}
The sampling theory approach being discussed in this section is based 
on~\cite{KG2003,WS2005}. Let the PSWFs be expanded as
\begin{equation}
\phi(t)=\sum_{n\in\field{Z}}\beta_n\psi_n(t),
\end{equation}
where $(\beta_n)_{n\in\field{Z}}$ are the unknown coefficients and $\psi_n(t)$ denotes 
normalized translates of the sinc function given by
\begin{equation}
\psi_n(t) = \sqrt{\frac{\sigma}{\pi}}\sinc[\sigma(t-t_n)]
=\sqrt{\frac{\sigma}{\pi}}\frac{\sin(\sigma t-n\pi)}{(\sigma t-n\pi)}
\end{equation}
so that the orthonormality condition reads as 
$\langle\psi_m,\psi_n\rangle_{\field{R}}=\delta_{mn}$. The 
coefficients $\beta_n$ can be worked out using the orthonormality property
or by using direct sampling at $t_n={n\pi}/{\sigma}$ so that 
$\beta_n=\sqrt{\pi/\sigma}\phi(t_n)$ which
is just another way of stating the Shannon-Whittaker sampling
theorem for $\sigma$-bandlimited functions~\cite{P1987}.

Using the expansion 
\begin{equation}\label{eq:sinc-ker-sinc}
\frac{\sin[\sigma(t-s)]}{\pi(t-s)}=\sum_{n\in\field{Z}}\psi_n(t)\psi_n(s),
\end{equation}
of the sinc-kernel in~\eqref{eq:PSWFs_EVP1}, the eigenvalue problem in~\eqref{eq:PSWFs_EVP1} 
can be discretized as:
\begin{equation}
\lambda\sum_{m\in\field{Z}}\beta_m\psi_m(t)
=\sum_{m\in\field{Z}}\psi_m(t)\sum_{n\in\field{Z}}\beta_n\int_{-1}^{1}ds\,\psi_m(s)\psi_n(s)
\end{equation}
Equating the coefficient of $\psi_m(t)$ both sides, we have
\begin{equation}\label{eq:sinc_evp}
\lambda\beta_m
=\sum_{n\in\field{Z}}A_{mn}\beta_n,
\end{equation}
where $A_{mn} = \langle\psi_m,\psi_n\rangle_{\Omega}$. 
These entries can be stated in terms of the sine and cosine integrals 
as follows:
\begin{equation}
A_{mn} 
=\left\{
\begin{aligned}
&\frac{2\sigma}{\pi}\frac{\sin^2\sigma}{(m^2\pi^2-\sigma^2)}\\
&+\frac{1}{\pi}\left[\Si(2\sigma-2m\pi)+\Si(2\sigma+2m\pi)\right],\,m=n,\\
&\frac{(-1)^{m+n}}{\pi^2(m-n)}[\Cin(2\sigma-2m\pi)-\Cin(2\sigma-2n\pi)\\
&-\Cin(2\sigma+m\pi)+\Cin(2\sigma+2n\pi)],\, m\neq n.
\end{aligned}\right.
\end{equation}
Using Cauchy-Schwarz inequality, we have
\begin{equation}
|A_{mn}|^2\leq\frac{2\sigma^2}{\pi(m^2\pi^2-\sigma^2)}
\frac{2\sigma^2}{\pi(n^2\pi^2-\sigma^2)},\quad |m|,|n|>\frac{\sigma}{\pi}.
\end{equation}
The operator $A$ can be shown to be a \emph{self-adjoint}, \emph{positive definite} and 
\emph{compact} operator on $\ell^2$ so that it admits of an orthonormal sequence of 
eigenvectors which 
correspond to positive eigenvalues~\cite{WS2005}. Various properties of the
PSWFs can also be deduced from the fact that they are eigenfunctions of the
operator $A$ as was done in~\cite{KG2003,WS2005}. 

% We conclude this section with
% the remark that the linear system in~\eqref{eq:sinc_evp} can be split into two
% parts where one corresponds to the odd and the other to the even parity PSWFs 
% by introducing the odd and even parity basis functions as
% follows: For $n\in\field{Z}_+$, define
% \begin{equation}
% \begin{split}
% &\psi_n^{(+)}(t)=\frac{1}{\sqrt{2}}[\psi_n(t)+\psi_{-n}(t)],\\
% &\psi_n^{(-)}(t)=\frac{1}{\sqrt{2}}[\psi_n(t)-\psi_{-n}(t)],
% \end{split}
% \end{equation}
% with $\psi_0^{(+)}(t)=\psi_0(t)$. The sinc-kernel in~\eqref{eq:sinc-ker-sinc} can
% be rewritten as
% \begin{equation}
% \frac{\sin[\sigma(t-s)]}{\pi(t-s)}=
% \sum_{n\geq0}\psi^{(+)}_n(t)\psi^{(+)}_n(s)
% +\sum_{n>0}\psi^{(-)}_n(t)\psi^{(-)}_n(s).
% \end{equation}
% Observing that the odd basis functions do not ``interact'' with the even ones,
% it is straightforward to decouple the two systems. Rest of the details are
% entirely similar to what has been carried out above which can be applied
% individually to the odd and the even parity PSWFs. In our tests, the numerical
% conditioning improves as a result of this splitting.

\section{Extrapolation of Bandlimited Signals}
\label{sec:Xtrap-BL}

The extrapolation problem for bandlimited signals known on $\Omega=(-1,1)$ requires the 
solution of the Fredholm equation~\cite{Cadzow1979} 
\begin{equation}\label{eq:fredholm}
x(t) = \int_{\Omega}\frac{\sin{\sigma(t-s)}}{\pi(t-s)}y(s)ds,\quad t\in\Omega,
\end{equation}
where $x(t)$ is known in the interval $\Omega$ and $y(t)$ is an unknown signal. The 
extrapolation of $x(t)$ for $t\in\field{R}\setminus \Omega$ is carried 
out using the same equation once $y(t)$ is determined on $\Omega$.

The discretization of this equation can be accomplished by using the
representation in~\eqref{eq:sinc_sphj} so that
\begin{equation}
x(t) = \sum_{n=0}^{\infty}\ovl{j}_n(t)\int_{\Omega}\ovl{j}_n(s)y(s)ds,
\end{equation}
The overlap integrals involving $y(t)$ constitute the unknowns; therefore, let 
$\hat{y}_n=\langle y,\ovl{j}_n\rangle_{\Omega}$ and 
define $\hat{x}_n=\langle x,\ovl{j}_n\rangle_{\Omega}$ to obtain the following
linear system of equations
\begin{equation}\label{eq:cadzow_onestep}
\hat{x}_m = \sum_{n}J_{mn}\hat{y}_n.
\end{equation}
The solution of the discrete system can be obtained in terms of the eigenvectors
of the symmetric, positive definite matrix $J$. The eigenvectors of this matrix
are precisely the PSWFs discussed in the last section. In the 
discrete form, the solution can be stated by writing the expansion of 
$J^{-1}$ using its eigenvectors 
\begin{equation}\label{eq:cadzow_sol}
\hat{\vv{y}} =
\sum_n\left(\frac{\hat{\vv{x}}\cdot\vs{\beta}_{n}}{\lambda_n}\right)\vs{\beta}_{n},
\end{equation}
provided 
\begin{equation}\label{eq:picard}
\sum_n
\left|\frac{\hat{\vv{x}}\cdot\vs{\beta}_{n}}{\lambda_n}\right|^2 <\infty.
\end{equation}
This is the special case of Picard's theorem~\cite[Chap.~VI]{P1966}.
For $n>2\sigma/\pi$, the eigenvalues $\lambda_n$ show a sharp decrease
rendering the sum~\eqref{eq:picard} extremely sensitive to the errors in
the $|\hat{\vv{x}}\cdot\vs{\beta}_n|$ making this problem inherently ill-posed. 
Thus, in order to obtain any meaningful solution of the problem, one must turn to 
regularization techniques for solving the linear system 
in~\eqref{eq:cadzow_onestep}.
\begin{rem}
In addition to the fact that $J$ has no bounded inverse, the regularity property of 
$x(t)$ dictates the quality of approximation
\begin{equation}
x_N(t)=\sum_{n=0}^{N-1}\alpha_n \ovl{j}_n(t).
\end{equation}
When $x\in\fs{B}^2_{\sigma}$ is known for all $t\in\field{R}$, one chooses 
$\alpha_n=\langle x,\ovl{j}_n\rangle_{\field{R}}$ so that
$\|x\|_2=\|\vs{\alpha}\|_{\ell^2}<\infty$. Choosing $N$ such that 
$2N+2>e\sigma|t|$ and taking into account the inequality
$|\ovl{j}_n(t)|\leq (e\sqrt{\sigma|t|}/2)[e\sigma|t|/(2n+2)]^{n}$,
we have
\begin{equation}
|x(t)-x_N(t)|\leq
\frac{e\sqrt{\sigma|t|}\,\|x\|_2}{2\left[1-\left(\frac{e\sigma|t|}{2N+2}\right)^2\right]^{1/2}}
\left(\frac{e\sigma|t|}{2N+2}\right)^{N}.
\end{equation}
Let $X(\xi)$ denote that Fourier transform of $x(t)$ and $X_N(\xi)$ denote that Fourier transform 
of $x_N(t)$. If $X(\xi)\in\fs{C}^{p}(-\sigma,\sigma)$ with $p>N$, then estimates for the 
truncation error that holds uniformly for all $t\in\field{R}$ can be obtained as 
follows~\cite{P1976}: Starting with 
\begin{equation*}
X_N(\sigma\xi) = \sum_{n=0}^{N-1}\beta_nP_n(\xi),\quad\xi\in\Omega,
\end{equation*}
we have
\begin{equation*}
x(t)-x_N(t)=\frac{\sigma^{N+1}}{2\pi N!!}
\int_{-1}^{1}X^{(N)}(\sigma\zeta_{\xi})P_N(\xi)e^{i\sigma\xi t}d\xi
\end{equation*}
where $\zeta_{\xi}\in(-1,1)$ so that
\begin{equation}
|x(t)-x_N(t)|\leq\frac{\sigma^{N+1}\|X^{(N)}\|_{\infty}}{\pi N!!\sqrt{2N+1}},\quad\forall t\in\field{R}.
\end{equation}
When $N>p$, estimates can be similarly obtained that guarantee a uniform truncation error 
that scale as $\bigO{N^{-p}}$ (see Canuto~\et~\cite[Chap.~5]{CHQZ2007}).
\end{rem}

\subsection{Tikhonov regularization}\label{sec:Tikh-xtrap}
\label{sec:Xtrap-BL-regu}
The Tikhonov regularization~\cite{TA1977} for the Fredholm
equation~\eqref{eq:fredholm} in the discretized form~\cite{H1990,H1992} can be stated as the minimization
of the following function of $\hat{\vv{y}}=(\hat{y}_0,\hat{y}_1,\ldots)^{\tp}$:
\begin{equation}\label{eq:Tikh-standard}
H(\hat{\vv{y}}) =
\|J\hat{\vv{y}}-\hat{\vv{x}}\|^2_{\ell^2}
+\mu^2\|\hat{\vv{y}}\|^2_{\ell^2},
\end{equation}
where we have used the fact that $J$ is symmetric and 
$\hat{\vv{x}}=(\hat{x}_0,\hat{x}_1,\ldots)^{\tp}$ with 
$\hat{x}_n=\langle x,\ovl{j}_n\rangle_{\Omega}$. This is known as the 
\emph{standard form} of Tikhonov regularization. The 
minimization problem for $\hat{\vv{y}}$ boils down to the solution of the 
linear system given by~\cite[Chap.~5]{H1998}
\begin{equation}\label{eq:Euler}
\left(\mu^2 I + J^2\right)\hat{\vv{y}}=J\hat{\vv{x}}.
\end{equation}
The solution can be stated in terms of the eigenvectors of $J$ as
\begin{equation}\label{eq:cadzow_sol_reg1}
\hat{\vv{y}}(\mu) = 
\sum_n\frac{\lambda_n}{\lambda^2_n+\mu^2}(\hat{\vv{x}}\cdot\vs{\beta}_{n})\vs{\beta}_{n}.
\end{equation}
In order to specify a sufficient criteria for the regularized solution to be
`close' to the actual solution under the Picard's condition~\eqref{eq:picard},
consider
\begin{equation}
\hat{\vv{y}}_{\text{exact}}-\hat{\vv{y}}(\mu)
=\sum_n\frac{\mu^2}{(\lambda^2_n+\mu^2)\lambda_n}(\hat{\vv{x}}\cdot\vs{\beta}_{n})\vs{\beta}_{n}.
\end{equation}
The right hand side of this equation is bounded on account of the Picard's
condition. For any $\epsilon>0$, we have
\begin{equation*}
\begin{split}
\|\hat{\vv{y}}_{\text{exact}}-\hat{\vv{y}}(\mu)\|^2_{\ell^2}
&\leq\sum_n\frac{\mu^{\epsilon}}{(\lambda^2_n+\mu^2)^{\epsilon/2}}
\left(\frac{\mu^2}{\lambda^2_n+\mu^2}\right)^{2-\epsilon/2}
\left|\frac{(\hat{\vv{x}}\cdot\vs{\beta}_{n})}{\lambda_n}\right|^2\\
&\leq\mu^{\epsilon}\sum_n\frac{|\hat{\vv{x}}\cdot\vs{\beta}_{n}|^2}{\lambda^{2+\epsilon}_n}.
\end{split}
\end{equation*}
Therefore, for the class of signals which satisfy the estimate 
\begin{equation}
\sum_n\frac{|\hat{\vv{x}}\cdot\vs{\beta}_{n}|^2}{\lambda^{2+\epsilon}_n}<\infty
\end{equation}
for some $\epsilon>0$, the convergence behavior with respect to the regularization
parameter can be explicitly determined to be
\begin{equation}
\|\hat{\vv{y}}_{\text{exact}}-\hat{\vv{y}}(\mu)\|_{\ell^2}=\bigO{\mu^{\epsilon/2}}.
\end{equation}

Truncation to a finite dimensional system can be achieved using the $N\times\infty$ projection operator 
$\mathcal{P}_N=(I_{N},0_{N\times\infty})$ where $I_N$ is the $N\times N$
identity matrix. Let 
${J}_N=\mathcal{P}_N J \mathcal{P}^{\tp}_N$ and define $\hat{\vv{y}}_N$ such that
\begin{equation}
\left[\mu^2 I_N + J^2_N\right]\hat{\vv{y}}_N
=J_N\mathcal{P}_N\hat{\vv{x}}.
\end{equation}
In the following, we suppress the dependence of the regularized solution $\hat{\vv{y}}$ and its
finite dimensional approximation $\hat{\vv{y}}_{N}$ on the regularization
parameter for the sake of brevity of presentation. Consider
\begin{multline}
\left({\mu}^{2} I + J^2\right)(\hat{\vv{y}}-\mathcal{P}^{\tp}_N\hat{\vv{y}}_N)+
\left[J^2 - \mathcal{P}^{\tp}_NJ^2_N\mathcal{P}_N\right]\mathcal{P}^{\tp}_N\hat{\vv{y}}_N\\
=\left(J-\mathcal{P}_N^{\tp}J_N\mathcal{P}_N\right)\hat{\vv{x}}
+\mathcal{P}^{\tp}_NJ_N\mathcal{P}_N\left(\hat{\vv{x}}-\mathcal{P}_N^{\tp}\mathcal{P}_N\hat{\vv{x}}\right),
\end{multline}
where $\mathcal{P}_N^{\tp}\mathcal{P}_N\hat{\vv{x}}$ agrees with $\hat{\vv{x}}$
up to first $N$ coefficients and has zero entries thereafter, and,
$\mathcal{P}_N^{\tp}J_N\mathcal{P}_N$ agrees with $J$ up to the first $N\times
N$ block and has zero entries thereafter. Our goal is to
obtain the error estimate 
$\|\hat{\vv{y}}-\mathcal{P}^{\tp}_N\hat{\vv{y}}_N\|_{\ell^2}$ in terms of
$\|J-\mathcal{P}_N^{\tp}J_N\mathcal{P}_N\|_s$ and
$\|\hat{\vv{x}}-\mathcal{P}_N^{\tp}\mathcal{P}_N\hat{\vv{x}}\|_{\ell^2}$. From
$\|J\|_s<1$ and $\|\mathcal{P}_N^{\tp}J_N\mathcal{P}_N\|_s<1$, we have
$\|J^2-\mathcal{P}_N^{\tp}J^2_N\mathcal{P}_N\|_s\leq2\|J-\mathcal{P}_N^{\tp}J_N\mathcal{P}_N\|_s$
which yields
\begin{multline}\label{eq:finite-d-error}
\|\hat{\vv{y}}-\mathcal{P}_N^{\tp}\hat{\vv{y}}_N\|_{\ell^2}
\leq\left\|(\mu^2 I+J^2)^{-1}\right\|_s\\
\times[(\|\hat{\vv{x}}\|_{\ell^2}+2\|\hat{\vv{y}}_N\|_{\ell^2})
\|J-\mathcal{P}_N^{\tp}J_N\mathcal{P}_N\|_s\\
+\|\hat{\vv{x}}-\mathcal{P}^{\tp}_N\mathcal{P}_N\hat{\vv{x}}\|_{\ell^2}].
\end{multline}
The estimates for  
$\|J-\mathcal{P}_N^{\tp}J_N\mathcal{P}_N\|_s$ and 
$\|\hat{\vv{x}}-\mathcal{P}^{\tp}_N\mathcal{P}_N\hat{\vv{x}}\|_{\ell^2}$
are provided in the following proposition:
\begin{prop}
For any positive integer $N$, define
\begin{equation}
\Gamma_N=\left(\frac{e\sigma}{2}\right)
\frac{\zeta^{2N}_N}{\left[1-\zeta_N^2\right]},\quad
\zeta_N=\frac{e\sigma}{(2N+2)}.
\end{equation}
If $N$ is chosen such that $\zeta_N<1$, the estimates
\begin{equation}\label{eq:x-xN}
\begin{split}
\|J-\mathcal{P}_N^{\tp}J_N\mathcal{P}_N\|_s\leq\Gamma_N,\quad 
\|\hat{\vv{x}}-\mathcal{P}^{\tp}_N\mathcal{P}_N\hat{\vv{x}}\|_{\ell^2}
\leq\|x\|_{\fs{L}^2(\Omega)}\Gamma_N^{1/2},
\end{split}
\end{equation}
hold.
\end{prop}
\begin{proof}
Proceeding with the Frobenius norm and using the estimate in~\eqref{eq:jmn-estimate2}, we have
\begin{equation}
\|J-\mathcal{P}_N^{\tp}J_N\mathcal{P}_N\|^2_{\text{F}}
\leq\left(\frac{e\sigma}{2}\right)^2
\left[\sum_{n\geq N}^{\infty}\left(\frac{e\sigma}{2n+2}\right)^{2n}\right]^2.
\end{equation}
Replacing $n$ by $N$ and summing the resulting geometric series yields the first estimate. In 
order to prove the second estimate, we first obtain an estimate for $\hat{x}_n$ as follows:
\begin{equation}
|\hat{x}_n|\leq\|x\|_{\fs{L}^2(\Omega)}\|\ovl{j}_n\|_{\fs{L}^2(\Omega)}
\leq\|x\|_{\fs{L}^2(\Omega)}\sqrt{\frac{e\sigma}{2}}\left(\frac{e\sigma}{2n+2}\right)^{n}.
\end{equation}
Following as in the first case, we have
\begin{equation*}
\begin{split}
\frac{\|\hat{\vv{x}}-\mathcal{P}^{\tp}_N\mathcal{P}_N\hat{\vv{x}}\|_{\ell^2}}{\|x\|_{\fs{L}^2(\Omega)}}
&\leq\sqrt{\frac{e\sigma}{2}}
\left[\sum_{n\geq N}\left(\frac{e\sigma}{2n+2}\right)^{2n}\right]^{1/2}
\leq\Gamma_N^{1/2}.
\end{split}
\end{equation*}
\end{proof}
From the proof above, we also conclude
\begin{equation}
\|\hat{\vv{x}}\|_{\ell^2}\leq\|x\|_{\fs{L}^2(\Omega)}
\sqrt{\frac{e\sigma}{2}}\left[\max_{0\leq n\leq N}\left(\frac{e\sigma}{2n+2}\right)^{n}
+\frac{1}{\sqrt{1-\zeta_N^2}}\right].
\end{equation}
Observing 
\begin{equation*}
\|\hat{\vv{y}}_N\|_{\ell^2}
\leq\left\|(\mu^2 I_N+J_N^2)^{-1}\right\|_s\|\hat{\vv{x}}_N\|_{\ell^2}
\leq\left\|(\mu^2 I+J^2)^{-1}\right\|_s\|\hat{\vv{x}}\|_{\ell^2},
\end{equation*}
it follows from~\eqref{eq:finite-d-error} that
\begin{equation}
\|\hat{\vv{y}}-\mathcal{P}_N^{\tp}\hat{\vv{y}}_N\|_{\ell^2}=\bigO{\zeta_N^{2N}/\mu^2}+\bigO{\zeta_N^{N}/\mu}
=\bigO{\zeta_N^{N}/\mu}.
\end{equation}
Therefore, if $\mu=\mu_N$ varies with $N$ in such a way that the ratio
$(\zeta_N/\mu_N)$
remains bounded, we have
\begin{multline}
\|\hat{\vv{y}}_{\text{exact}}-\mathcal{P}_N^{\tp}\hat{\vv{y}}_N\|_{\ell^2}
=\bigO{\mu_N^{\epsilon/2}}+\bigO{\frac{\zeta_N^{N}}{\mu_N}}\\
=\bigO{\mu_N^{\epsilon/2}}+\bigO{\zeta_N^{N}}.
\end{multline}
The advantage of the spherical Bessel basis functions is apparent from the fact
that the contribution of the discretization error to the reconstruction error is
exponentially small. Let us also note that the error constants which are missing
in the Landau's `O' notation increases asymptotically as some positive power of
$\sigma$. 

Given that the bandlimited signals are smooth functions, we may 
introduce extra ``penalty'' for the lack of smoothness by considering 
Sobolev norm type penalty terms in the minimization problem. To this end, 
let us show that it is possible to construct a differentiation matrix 
which acts on the discrete representation $\vs{\beta}\in\ell^2$ of any 
signal $x(t)\in\fs{B}_{\sigma}^2$ to compute $x'(t)\in\fs{B}_{\sigma}^2$ 
which we note satisfies the estimate 
$\|x'\|_2\leq \sigma\|x\|_2$~\cite[Chap.~2]{M2001}. The first 
order derivative of the spherical Bessel functions obey the recurrence relation
\begin{equation}
(2n+1)j'_{n}(\sigma t) = nj_{n-1}(\sigma t)-(n+1)j_{n+1}(\sigma t),\quad n>0,
\end{equation}
with $j'_0(\sigma t)= - j_1(\sigma t)$ so that
\begin{equation}\label{eq:diff_psf_bessel}
x'(t)=\sum_{n=0}^{\infty}\beta_n \ovl{j}'_n(t)
=\sum_{n=0}^{\infty}\beta'_n\ovl{j}_n(t),
\end{equation}
where $\beta'_n$ are to be determined. Equating the coefficient 
of $\ovl{j}_m( t)$, we have
\begin{equation}\label{eq:dbeta}
\beta'_m 
=\frac{\sigma (m+1)\beta_{m+1}}{\sqrt{(2m+1)(2m+3)}}
-\frac{\sigma m\beta_{m-1}}{\sqrt{(2m-1)(2m+1)}},
\end{equation}
which yields the first order differentiation matrix $\sigma D^{(1)}$ that turns 
out to be a skew-symmetric banded matrix with its non-zero entries given by
\begin{equation}\label{eq:D1}
%\begin{split}
D^{(1)}_{m,m+1}=\frac{m+1}{\sqrt{(2m+1)(2m+3)}}.
%&D^{(1)}_{m,m-1}=-\frac{m}{\sqrt{(2m-1)(2m+1)}}.
%\end{split}
\end{equation}
Let us state the following result for discrete representation of the 
the first derivative:

\begin{figure*}[th!]
\centering
\includegraphics[scale=1]{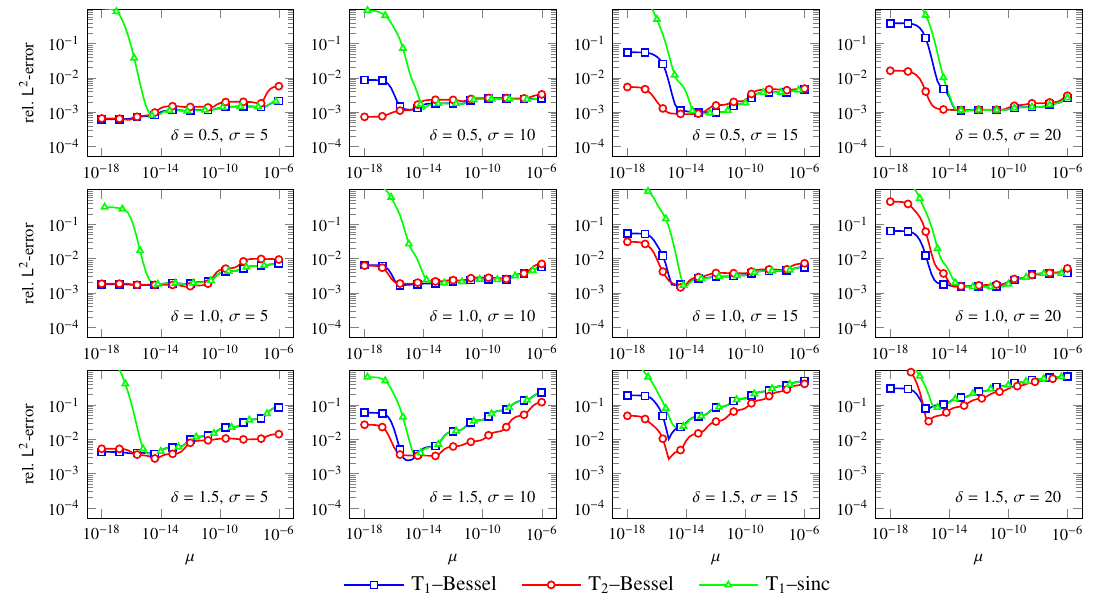}
\caption{\label{fig:xtrap}The figure shows the error 
(quantified by~\eqref{eq:error}) as a function of the
regularization parameter ($\mu$) for various bandlimited
extrapolation algorithms for the signal $x(t)$ as defined 
in~\eqref{eq:xI} with $p=1$. The signal is assumed to be 
known in $\Omega$ and the
extrapolation error is computed over $\Omega'=(-5,5)$.}
\end{figure*}

\begin{figure*}[th!]
\centering
\includegraphics[scale=1]{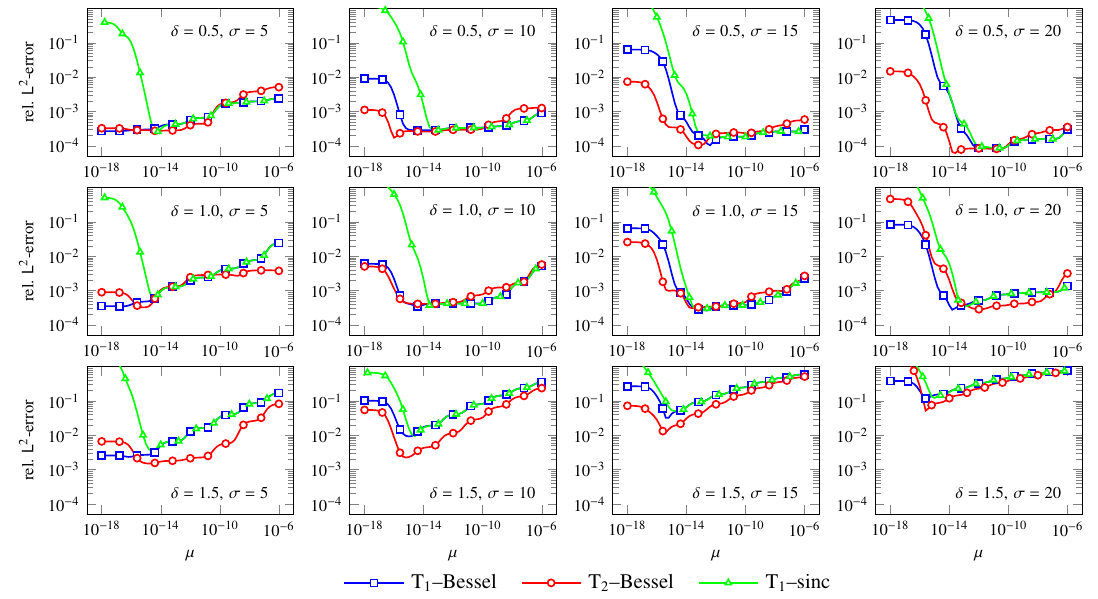}
\caption{\label{fig:xtrap2}The figure shows the error 
(quantified by~\eqref{eq:error}) as a function of the
regularization parameter ($\mu$) for various bandlimited
extrapolation algorithms for the signal $x(t)$ as defined 
in~\eqref{eq:xI} with $p=2$. The signal is assumed to be 
known in $\Omega$ and the
extrapolation error is computed over $\Omega'=(-5,5)$.}
\end{figure*}

\begin{figure*}[th!]
\centering
\includegraphics[scale=1]{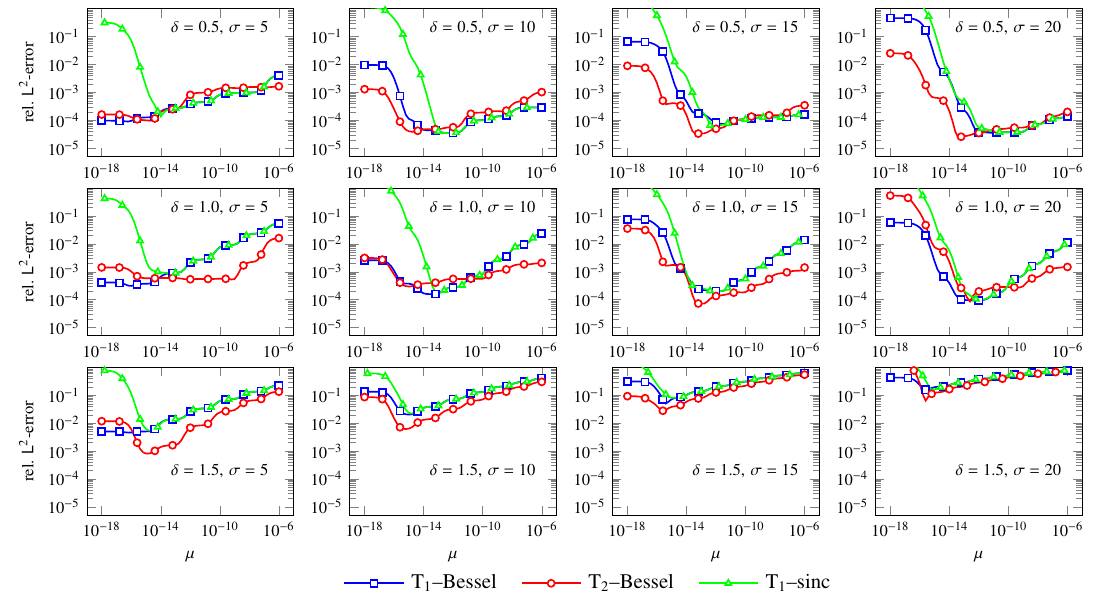}
\caption{\label{fig:xtrap3}The figure shows the error 
(quantified by~\eqref{eq:error}) as a function of the
regularization parameter ($\mu$) for various bandlimited
extrapolation algorithms for the signal $x(t)$ as defined 
in~\eqref{eq:xI} with $p=3$. The signal is assumed to be 
known in $\Omega$ and the
extrapolation error is computed over $\Omega'=(-5,5)$.}
\end{figure*}
\begin{lemma}
Let the representation of any signal $x\in\fs{B}^2_{\sigma}$ in 
the normalized spherical Bessel function basis be $\vs{\beta}\in\ell^2$ 
and that of its derivative $x'$ (which also belongs 
to $\fs{B}^2_{\sigma}$) be $\vs{\beta}'\in\ell^2$; then, 
we have $\vs{\beta}'=\sigma D^{(1)}\vs{\beta}$ such that 
$\|\vs{\beta}'\|_{\ell^2}\leq\sigma\|\vs{\beta}\|_{\ell^2}$.
\end{lemma}
\begin{proof}
The relation $\vs{\beta}'=\sigma D^{(1)}\vs{\beta}$ follows from 
the recurrence relation~\eqref{eq:dbeta}. Using the Cauchy-Schwarz 
inequality in the aforementioned recurrence relation, we have
$|\beta'_m|^2\leq\sigma^2(|\beta_{m+1}|^2+|\beta_{m-1}|^2)(C_{m+1}+C_{m})$
where $C_m=m^2/(4m^2-1)$. Observing $C_m=m^2/(4m^2-1)\leq 1/4$ for 
all $m\geq0$, we have 
$2|\beta'_m|^2\leq\sigma^2(|\beta_{m+1}|^2+|\beta_{m-1}|^2)$ 
which yields $\|\vs{\beta}'\|_{\ell^2}\leq\sigma\|\vs{\beta}\|_{\ell^2}$.
\end{proof}
The lemma shows that $\|D^{(1)}\|_s\leq1$. Further, using the recurrence 
relation~\eqref{eq:dbeta}, it is also possible to 
construct the second order differentiation matrix $\sigma^2D^{(2)}$ which 
turns out to be a symmetric banded matrix with diagonal entries given by 
\begin{equation}
D^{(2)}_{m,m}=-\frac{2m^2+2m-1}{(2m-1)(2m+3)},
\end{equation}
and the nonzero off-diagonal elements given by
\begin{equation}
D^{(2)}_{m,m+2}=\frac{(m+1)(m+2)}{(2m+3)\sqrt{(2m+1)(2m+5)}}.
\end{equation}
Now, we can consider
a minimization problem of the form
\begin{equation}
H(\hat{\vv{y}}) =
\|J\hat{\vv{y}}-\hat{\vv{x}}\|^2_{\ell^2}
+\mu^2\left(\|\hat{\vv{y}}\|^2_{\ell^2}
+\|D^{(1)}\hat{\vv{y}}\|^2_{\ell^2}
+\|D^{(2)}\hat{\vv{y}}\|^2_{\ell^2}
\right).
\end{equation}
The quantity in parenthesis above is similar to a 
Sobolev norm. As before, this minimization problem is equivalent to 
solving the linear system given by  
\begin{equation}
\left[J^2+\mu^2\left(I-{D}^{(1)}{D}^{(1)}+{D}^{(2)}{D}^{(2)}\right)\right]\hat{\vv{y}}
=J\hat{\vv{x}},
\end{equation}
where we used the skew-symmetric nature of $D^{(1)}$. This equation can be solved 
using the generalized SVD (see~\cite{H1992,H1998}). 

Finally, let us remark that the Tikhonov regularization can also be discussed in
the discrete framework based on the translates of the sinc function. These details are 
being omitted here because the line of reasoning is entirely similar. However, with regard to 
the translates of the sinc function, let us note that it is considerably 
harder to implement the Sobolev norm in this basis; therefore, we restrict 
ourselves to the standard form of the Tikhonov regularization, 
i.e.~\eqref{eq:Tikh-standard}, in this case.

%%%%%%%%%%%%%%%%%%%%%%%%%%%%%%%%%%%%%%%%%%%%%%%%%%%%%%%%%%%%%%%%%%%%%%%%%%%%%%%%
%%%%%%%%%%%%%%%%%%%%%%%%%%%%%%%%%%%%%%%%%%%%%%%%%%%%%%%%%%%%%%%%%%%%%%%%%%%%%%%%
\subsection{Numerical examples}\label{sec:xtrap-result}
For the purpose of numerical tests, we first label the algorithms being tested
as follows:
\begin{itemize}
\item T$_1$--Bessel: Standard Tikhonov regularization with spherical Bessel
functions as basis functions.
\item T$_2$--Bessel: Tikhonov regularization where regularity is enforced via a
Sobolev norm with spherical Bessel functions as basis functions. 
\item T$_1$--sinc: Standard Tikhonov regularization with translates of sinc
function as basis functions.
\end{itemize}

\begin{figure*}[th!]
\centering
\includegraphics[scale=1]{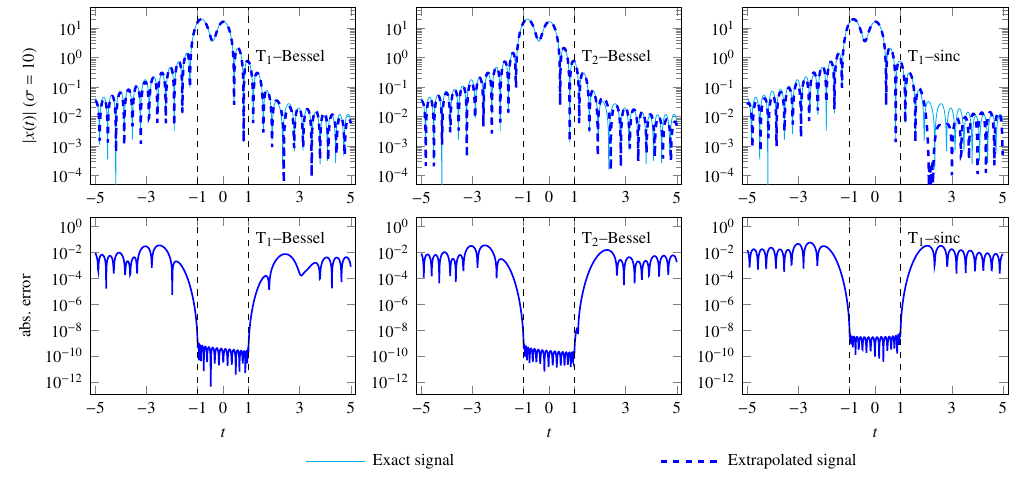}
\caption{\label{fig:xtrap-sig} The figure shows a comparison of the
extrapolated and the exact signal~\eqref{eq:xI} for $p=1$, $\sigma=10$ and $\delta=1.0$ 
using different algorithms. The point-wise error, 
$|x(t)-x_{\text{extrap.}}(t)|$, is plotted in the second row. The regularization 
parameter chosen to be $\mu=10^{-15}$ for T$_1$- and T$_2$-Bessel while $\mu=10^{-14}$ for T$_1$-sinc.}

\end{figure*}
\begin{figure*}[th!]
\centering
\includegraphics[scale=1]{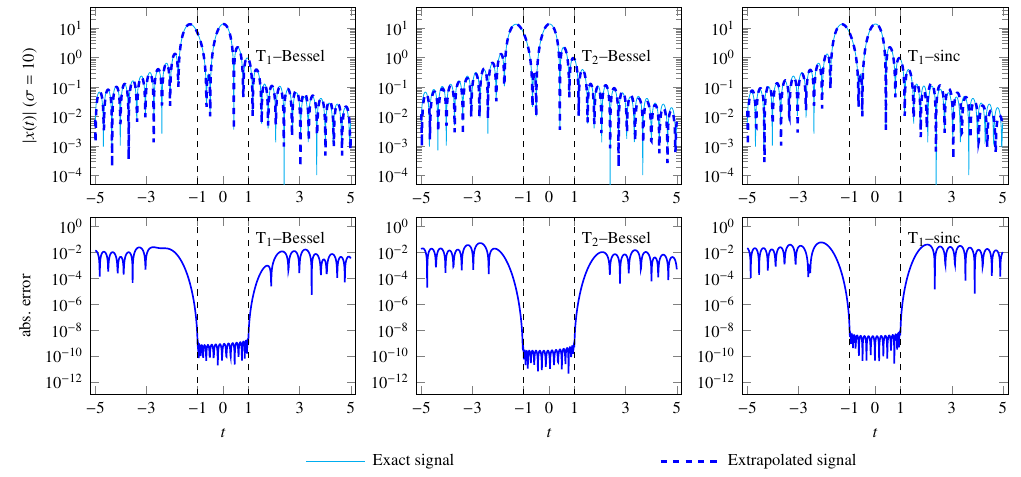}
\caption{\label{fig:xtrap-sig2} The figure shows a comparison of the
extrapolated and the exact signal~\eqref{eq:xI} for $p=1$, $\sigma=10$ and $\delta=1.5$ 
using different algorithms. The point-wise error, 
$|x(t)-x_{\text{extrap.}}(t)|$, is plotted in the second row. The regularization 
parameter chosen to be $\mu=10^{-15}$ for T$_1$- and T$_2$-Bessel while $\mu=10^{-14}$ for T$_1$-sinc.}
\end{figure*}
The number of basis functions taken is $N=400$. The number of Legendre--Gauss--Lobatto 
quadrature nodes is chosen to be
$N_{\text{quad.}}=1600$. The bandlimiting parameter $\sigma$ is chosen from
$\{10,20\}$. The regularization parameter $\mu$ is allowed to vary in a range
for which we look at the error in the extrapolated signal. The error is
quantified by a relative $\fs{L}^2(\Omega')$-norm where $\Omega'=(-5,5)$:
\begin{equation}\label{eq:error}
e_{\text{rel.}}={\|x_{\text{exact}}-x_{\text{extrap.}}\|_{\fs{L}^2(\Omega')}}
/{\|x_{\text{exact}}\|_{\fs{L}^2(\Omega')}},
\end{equation}
where $x_{\text{exact}}$ denotes the exact signal and $x_{\text{extrap.}}$
denotes the extrapolated signal. Note that the signal is assumed to be known 
in $\Omega$. In this paper, we do not
present a study of the effectiveness of different methods of finding the optimal
regularization parameter $\mu$; however, let us mention that the L-curve method
seems to perform satisfactorily. Finally let us note that the discrete problem 
can also be solved using the \emph{mean gradient descent} approach recently proposed 
by Rajora~\et~\cite{RBK2019} which tends to find a solution of the ill-posed problem by balancing 
the residue term and the penalty term.

\begin{figure*}[th!]
\centering
\includegraphics[scale=1]{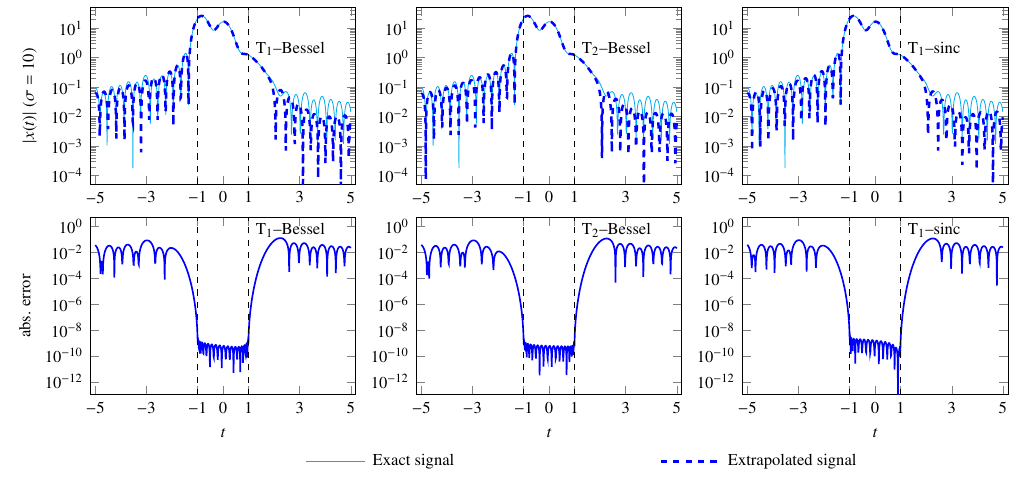}
\caption{\label{fig:xtrap-sig3} The figure shows a comparison of the
extrapolated and the exact signal~\eqref{eq:xII} for $p=1$, $\sigma=10$ and $\delta=1.0$ 
using different algorithms. The point-wise error, 
$|x(t)-x_{\text{extrap.}}(t)|$, is plotted in the second row. The regularization 
parameter chosen to be $\mu=10^{-15}$ for T$_1$- and T$_2$-Bessel while $\mu=10^{-14}$ for T$_1$-sinc.}
\end{figure*}

\begin{figure*}[th!]
\centering
\includegraphics[scale=1]{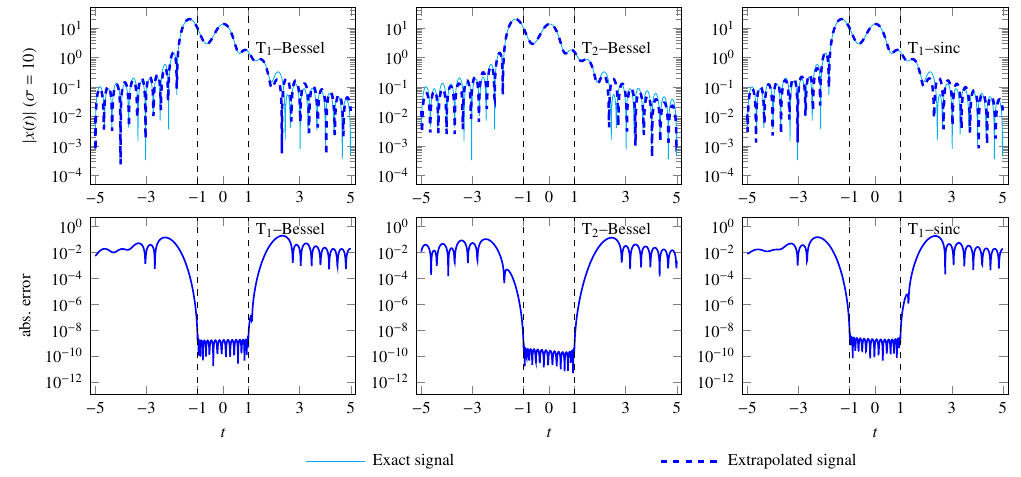}
\caption{\label{fig:xtrap-sig4} The figure shows a comparison of the
extrapolated and the exact signal~\eqref{eq:xII} for $p=1$, $\sigma=10$ and $\delta=1.5$ 
using different algorithms. The point-wise error, 
$|x(t)-x_{\text{extrap.}}(t)|$, is plotted in the second row. The regularization 
parameter chosen to be $\mu=10^{-15}$ for T$_1$- and T$_2$-Bessel while $\mu=10^{-14}$ for T$_1$-sinc.}
\end{figure*}

For the purpose of testing, we consider signals composed out of products of 
the sinc and the Bessel functions (denoted by $J_p(\cdot)$ where $p$ is the order):
\begin{equation}\label{eq:xI}
x(t)=\sum_{l=0}^{n-1}\sinc(\sigma\kappa_lt)\frac{J_{p}\left(\sigma(1-\kappa_l)(t-\tau_l)\right)}
{(t-\tau_l)^{p}},
\end{equation}
where $\tau_l\in[-\delta,\delta]\, (\delta>0)$ and $\kappa_l\in[0,1]$ are chosen 
randomly from a uniform distribution. Evidently, the Fourier spectrum 
$X(\xi)\in\fs{C}^{p}(-\sigma,\sigma)$ (which follows from the observation that 
$t^{p}x(t)\in\fs{L}^1(\field{R})$). Here $\delta$ 
controls whether the dominant part of the signal lies within $\Omega$. The 
parameters of the test cases are chosen as follows: $n=20$, $\delta\in\{0.5,1.0,1.5\}$, 
$\sigma\in\{5,10,15,20\}$ and $p\in\{1,2,3\}$. The results of the numerical experiments 
where we have plotted the error defined by~\eqref{eq:error} as a function of the 
regularization parameter are shown in Fig.~\ref{fig:xtrap}--\ref{fig:xtrap3} which 
correspond to different values of $p$. Based on the results it is possible to identify 
that allowing $\delta$ to be larger than $1$ (which corresponds to a situation where 
the dominant part of the signal falls outside $\Omega$) makes the extrapolation more 
challenging. The minimum achievable error in this regime (i.e. for $\delta>1$) tends 
to increase with increasing bandwidth ($2\sigma$). Also, the regularization scheme 
based on Sobolev norm type penalty tends to perform better than other regularization 
schemes in this regime. Finally, let us note that the discrete system based on the spherical 
Bessel functions does not give large error when the regularization parameter is 
arbitrarily chosen to be a quantity close to the machine epsilon. The quality of 
the extrapolated signal for the cases $\delta=1.0$ and $\delta=1.5$ can be assessed from 
Fig.~\ref{fig:xtrap-sig} and~\ref{fig:xtrap-sig2}, respectively. We conclude this discussion 
with a second closely related example where the sinc function is replaced by $\Si(t)/t$: 
\begin{equation}\label{eq:xII}
x(t)=\sum_{l=0}^{n-1}
\left[\frac{\Si(\sigma\kappa_lt)}{\sigma\kappa_lt}\right]
\frac{J_{p}\left(\sigma(1-\kappa_l)(t-\tau_l)\right)}
{(t-\tau_l)^{p}}.
\end{equation}
The extrapolated signal for the cases $\delta=1.0$ and $\delta=1.5$ are shown in 
Fig.~\ref{fig:xtrap-sig3} and~\ref{fig:xtrap-sig4}, respectively.

%%%%%%%%%%%%%%%%%%%%%%%%%%%%%%%%%%%%%%%%%%%%%%%%%%%%%%
%%%%%%%%%%%%%%%%%%%%%%%%%%%%%%%%%%%%%%%%%%%%%%%%%%%%%%%%%%%%%%%%%%%%%%%
%%%%%%%%%%%%%%%%%%%%%%%%%%%%%%%%%%%%%%%%%%%%%%%%%%%%%%%%%%%%%%%%%%%%%%%
%%%%%%%%%%%%%%%%%%%%%%%%%%%%%%%%%%%%%%%%%%%%%%%%%%%%%%%%%%%%%%%%%%%%%%%
%%%%%%%%%%%%%%%%%%%%%%%%%%%%%%%%%%%%%%%%%%%%%%%%%%%%%%%%%%%%%%%%%%%%%%%
%%%%%%%%%%%%%%%%%%%%%%%%%%%%%%%%%%%%%%%%%%%%%%%%%%%%%%%%%%%%%%%%%%%%%%%
%%%%%%%%%%%%%%%%%%%%%%%%%%%%%%%%%%%%%%%%%%%%%%%%%%%%%%%%%%%%%%%%%%%%%%%
\section{Conclusion}\label{sec:conclusion}
To conclude, we have devised a new degenerate kernel method for solving 
Fredholm integral equation associated with the bandlimited extrapolation
problem. The regularization procedure was applied at the discrete level. The choice of
using spherical Bessel functions as the basis functions improves the performance
of the numerical algorithms when compared to a similar method using the
translates of the sinc function. Further, it is simpler to compute the
derivative of the signal in the representation using spherical Bessel functions 
so that Sobolev norm type penalty terms can be employed in the regularization procedure. 
% trigger a \newpage just before the given reference
% number - used to balance the columns on the last page
% adjust value as needed - may need to be readjusted if
% the document is modified later
\IEEEtriggeratref{26}
% The "triggered" command can be changed if desired:
\bibliographystyle{IEEEtran}
%\bibliography{PSWF}
%\IEEEtriggercmd{\enlargethispage{-10cm}}
% Generated by IEEEtran.bst, version: 1.14 (2015/08/26)
\providecommand{\noopsort}[1]{}\providecommand{\singleletter}[1]{#1}%

\end{document}